\documentclass[a4paper,11pt]{amsart}
\usepackage{amsmath,amssymb,amsfonts,amsthm}

\usepackage[alphabetic]{amsrefs}
\usepackage{bm}
\usepackage{graphicx}
\usepackage{ascmac}
\usepackage[all]{xy}
\usepackage{amsthm}
\usepackage[top=25mm,bottom=25mm,left=30mm,right=30mm]{geometry}
\usepackage{tikz}
\usetikzlibrary{patterns}
\usepackage{pxpgfmark}
\usepackage{multirow}
\usetikzlibrary{intersections,arrows,calc,decorations.markings}
\newtheorem{theorem}{Theorem}[section]
\newtheorem{proposition}[theorem]{Proposition}

\newtheorem{lemma}[theorem]{Lemma}

\theoremstyle{definition}
\newtheorem{remark}[theorem]{Remark}

\newtheorem{definition}[theorem]{Definition}

\makeatletter
 
 \@addtoreset{equation}{section}
\makeatother

\def\cc{\mathbf{c}}

\def\ee{\mathbf{e}}

\title[Non-Leaving-Face property of exchange graphs]{Exchange graphs of cluster algebras have the non-leaving-face property}
\author{Changjian Fu}
\address{Changjian Fu\\Department of Mathematics\\SiChuan University\\610064 Chengdu\\P.R.China}
\email{changjianfu@scu.edu.cn}
\date{\today}
\author{Shengfei Geng}
\address{Shengfei Geng\\Department of Mathematics\\SiChuan University\\610064 Chengdu\\P.R.China}
\email{genshengfei@scu.edu.cn}
\author{Pin Liu}
\address{Pin Liu\\Department of Mathematics\\
	Southwest Jiaotong University\\
	610031 Chengdu \\
	P.R.China}
\email{pinliu@swjtu.edu.cn}
\dedicatory{Dedicated to Professor Jie Xiao on the Occasion of his 60th Birthday}

\begin{document}
\maketitle

\begin{abstract}
    The claim in the title is proved.
\end{abstract}

\section{Introduction}
 Generalized associahedra were introduced by Fomin and
 Zelevinsky in connection to cluster algebras of finite type. Ceballos and Pilaud \cite{CP2016} showed that all type $B$-$C$-$D$ associahedra have the non-leaving-face property, namely, any geodesic connecting two vertices in the graph of the polytope stays in the minimal face containing both. In fact, this property was already proven by Sleator, Tarjan and Thurston \cite{STT1988} for associahedra of type $A$, before the name was coined in \cite{CP2016}.  For a finite Coxeter system $(W,S)$, Williams \cite{Will2017} established the non-leaving-face property for $W$-permuation and $W$-associahedra. However, it is known that there are examples related to the associahedron which do not satisfy the non-leaving-face property (cf. \cite{CP2016}).
 
An important combinatorial invariant of a cluster algebra is its exchange graph. The abstract exchange graph was introduced in \cite{BY2013}, which unifies various exchange graphs arising from representation theory of finite-dimensional algebras, marked surfaces, cluster algebras and so on. An abstract exchange graph has the structure of a generalized polytope and then the non-leaving-face property can be formulated in this general situation. In \cite{BZ2019}, the non-leaving-face property for exchange graphs arising from unpunctured marked surfaces was established. 
 
 The main contribution of this note is Theorem \ref{t:cluster}, in which we prove that for any cluster algebra, the exchange graph has the non-leaving-face property. 

\section{Background on cluster algebras}

Fix a positive integer $n$. Let $\mathcal{F}$ be the field of  rational functions in $n$ indeterminates with coefficients in $\mathbb{Q}$. 

First of all, a \emph{labeled seed} is a pair $(\mathbf{x}, B)$, where
\begin{itemize}
	\item $\mathbf{x}=(x_1, \dots, x_n)$ is an $n$-tuple of elements of $\mathcal F$ forming a free generating set of $\mathcal F$, that is, $x_1, \cdots, x_n$ are algebraically independent and $\mathcal{F}=\mathbb{Q}(x_1,\cdots,x_n)$.
	\item $B=(b_{ij})$ is an $n \times n$ integer matrix which is \emph{skew-symmetrizable}, that is, there exists a positive integer diagonal matrix $S$ such that $SB$ is skew-symmetric. 
\end{itemize}
We have the terminology:
\begin{itemize}
\item $\mathbf{x}$ is the (labeled) \emph{cluster} of this seed;
\item $x_1,\dots,x_n$ are its cluster variables;
\item $B$ is the exchange matrix.
\end{itemize}

As in \cite{FZ2002}, let $k \in\{1,\dots, n\}$. The \emph{seed mutation $\mu_k$ in direction $k$} transforms  $(\mathbf{x}, B)$ into a new labeled seed $(\mathbf{x'}, B'):=\mu_k(\mathbf{x}, B)$ defined as follows:
\begin{itemize}
	\item The entries of $B'=(b'_{ij})$ are given by 
	\begin{align} \label{eq:matrix-mutation}
		b'_{ij}=\begin{cases}-b_{ij} &\text{if $i=k$ or $j=k$,} \\ 
			b_{ij}+\left[ b_{ik}\right] _{+}b_{kj}+b_{ik}\left[ -b_{kj}\right]_+ &\text{otherwise.}
		\end{cases}
	\end{align}where we use the notation $[b]_+=\max(b,0)$ for an integer $b$.
	\item The cluster variables $\mathbf{x'}=(x'_1, \dots, x'_n)$ are given by
	\begin{align}\label{eq:x-mutation}
		x'_j=\begin{cases}\dfrac{\mathop{\prod}\limits_{i=1}^{n} x_i^{[b_{ik}]_+}+\mathop{\prod}\limits_{i=1}^{n} x_i^{[-b_{ik}]_+}}{x_k} &\text{if $j=k$,}\\
			x_j &\text{otherwise.}
		\end{cases}
	\end{align}
\end{itemize}

Let $\mathbb{T}_n$ be the \emph{$n$-regular tree} whose edges are labeled by the numbers $1, \dots, n$ such that the $n$ edges emanating from each vertex have different labels. We write 
\begin{xy}(0,1)*+{t}="A",(10,1)*+{t'}="B",\ar@{-}^k"A";"B" \end{xy} 
to indicate that vertices $t,t'\in \mathbb{T}_n$ are joined by an edge labeled by $k$. 

A seed pattern is defined by assigning a labeled seed $(\mathbf{x}_t, B_t)$ to every vertex $t\in \mathbb{T}_n$, so that the seeds assigned to the end points of any edge \begin{xy}(0,1)*+{t}="A",(10,1)*+{t'}="B",\ar@{-}^k"A";"B" \end{xy} are obtained from each other by the seed mutation in direction $k$. Let
\[
\mathcal{X}_{\mathcal{A}}:=\bigcup_{t\in\mathbb{T}_n}\mathbf{x}_t
\]be the set of all cluster variables appearing in its seeds. The cluster algebra $\mathcal{A}$ is the $\mathbb{Q}$-subalgebra of $\mathcal{F}$ generated by $\mathcal{X}_{\mathcal{A}}$.

We now fix a vertex $s\in \mathbb{T}_n$ and a skew-symmetrizable integer $n \times n$  matrix $B=(b_{ij})=B_s$. To $s$ and $B$ we associate a family of integer vectors $\mathbf{c}_{j;t}^{B_s;s}\in\mathbb{Z}^n$ (c-vectors); here $j=1,\ldots,n$ and $t\in\mathbb{T}_n$. According to \cite{FZ2007}, the vectors $\mathbf{c}_{j;t}^{B_s;s}$ can be defined by the initial conditions
\[
\cc_{j;s}^{B_s;s}=\mathbf{e}_j\ (j=1,\ldots, n)
\]
together with the recurrence relations
\begin{align}\label{f:c-vector}
	\cc_{j;t'}^{B_s;s} =
	\begin{cases}
		-\cc_{j;t}^{B_s;s} & \text{if $j= k$;} \\[.05in]
		\cc_{j;t}^{B_s;s} + [b_{kj;t}]_+ \ \cc_{k;t}^{B_s;s} +b_{kj;t} [-\cc_{k;t}^{B_s;s}]_+
		& \text{if $j\neq k$}.
	\end{cases}
\end{align} for any \begin{xy}(0,1)*+{t}="A",(10,1)*+{t'}="B",\ar@{-}^k"A";"B" \end{xy} in $\mathbb{T}_n$. Here $\mathbf{e}_1,\cdots, \mathbf{e}_n$ are unit vectors in $\mathbb{Z}^n$ and for an integer matrix $C$, we write $[C]_+$ for the matrix
obtained from $C$ by applying the operation $c\mapsto [c]_+$  to all entries of $C$. The $n\times n$ integer matrix 
\[
C_t^{B_s;s}=(\mathbf{c}_{1;t}^{B_s;s},\dots,\mathbf{c}_{1;t}^{B_s;s})
\]is called the \emph{$C$-matrix} at vertex $t$ with respect to $s$.

A non-zero integer vector $\cc\in \mathbb{Z}^n$ is \emph{sign-coherent} if all entries are either non-negative or non-positive. The following statement was proven in \cite{GHKK}.  
\begin{theorem}[sign-coherence of $c$-vectors]
	For each vertex $t$ in $\mathbb{T}_n$, the matrix $C_t^{B_s;s}$ is inveritble over $\mathbb{Z}$ and each column vector is sign-coherent.
\end{theorem}

For a matrix $C\in \text{M}_n(\mathbb{Z})$ and an index $k$, we will denote the matrix obtained from $C$ by replacing all entries outside of the $k$-th row (resp. column) with zeros $C^{k\bullet}$ (resp. $C^{\bullet k}$). As noted in \cite{NZ}, the operations $C\mapsto [C]_+$ and $C\mapsto C^{k\bullet}$
commute with each other, making the notation $[C]_+^{k\bullet}$ (and $[C]_+^{\bullet k}$) unambiguous. So the sign-coherence of $c$-vectors implies that for every $k$, there is the sign $\epsilon_k(C_t^{B_s;s})=\pm 1$ such that $[-\epsilon_k(C_t^{B_s;s})C_t^{B_s;s}]_+^{\bullet k}=0$. For each $1\leq k\leq n$, 
let $J_k$ be the diagonal matrix obtained from the identity matrix by replacing the $(k,k)$-entry with $-1$. The following statement was proven in \cite[Proposition 1.4]{NZ}.
\begin{proposition}\label{p:trans-c-vector}
	Let $\xymatrix@C=0.5cm{s\ar@{-}[r]^k&s'}$ be an edge of $\mathbb{T}_n$. Then for any vertex $t\in\mathbb{T}_n$
	\[C_t^{B_{s'};s'}=(J_k+[-\epsilon_k(C_s^{-B_t;t})B_s]_+^{k\bullet})C_t^{B_s;s}.
	\]
	where $C_s^{-B_t;t}$ is the $C$-matrix at $s$ with respect to $t$ by assigning $-B_t$ to $t$.
\end{proposition}

\section{Non-Leaving-Face property of exchange graphs}

\subsection{Background on exchange graphs and non-leaving-face property}\label{ss:abstract-exchange-graph}
As in \cite{BY2013}, let $V$ be a non-empty set with a reflexive and symmetric relation $R$. Two elements $x$ and $y$ of $V$ are {\it compatible} if $(x,y)\in R$. Assume $(V,R)$ satisfies the following conditions:
\begin{itemize}
    \item All maximal subsets of pairwise compatible elements are finite and have the same cardinality, say $n$. We refer the subsets {\it clusters};
    \item Any subset of $n-1$ pairwise compatible elements is contained in precisely two clusters.
\end{itemize}

The {\it exchange graph} $\mathcal{G}$ of $(V,R)$ is defined to be the graph whose vertices are the clusters and where two vertices are joined by an edge if and only if their intersection has cardinality $n-1$.
A {\it face} $\mathcal{F}$ of $G$ is a full subgraph of $G$ such that
\begin{itemize}
    \item there is a subset $U$ of pairwise compatible elements;
    \item the vertices of $\mathcal{F}$ are precisely the clusters containing $U$ as a subset.
\end{itemize}
Clearly, $\mathcal{F}$ is uniquely determined by $U$ and we denote it by $\mathcal{F}_U$. The inclusion of sets induces a partial order on faces. Namely, we say $\mathcal{F}_U\leq \mathcal{F}_V$ if $V\subseteq U$.

For two cluster $v_1$ and $v_2$, we write for $\xymatrix@C=0.5cm{v_1\ar@{-}[r]&v_2}$ to indicate that they are linked by an edge. A path 
\[\xymatrix@C=0.5cm{v=v_1\ar@{-}[r]&v_2\ar@{-}[r]&\cdots\ar@{-}[r]&v_m=w}
\]
from $v$ to $w$ is a {\it geodesic} connecting vertices $v$ and $w$, if the length of the path is minimal. We remark that the graph $\mathcal{G}$ may not be connected.

We now recall non-leaving-face property from \cite{CP2016}.
\begin{definition}
If any geodesic connecting two vertices in $\mathcal{G}$ lies in the minimal face containing both, then we say that $\mathcal{G}$ has {\it non-leaving-face} property.
\end{definition}

The following definition is useful to investigate the non-leaving-face property (cf. \cite{STT1988,CP2016,Will2017,BZ2019}).

\begin{definition}\label{d:projection}
	Let $\mathcal{G}$ be an exchange graph and $\mathcal{F}$ a face of $\mathcal{G}$. We say a map $P:\mathcal{G}\to \mathcal{F}$ is a \emph{projection} if the following conditions satisfied
	\begin{enumerate}
		\item[(P1)] $P(v)$ is a vertex in $\mathcal{F}$ for each vertex $v\in \mathcal{G}$;
		\item[(P2)] $P(v)=v$ whenever $v\in \mathcal{F}$;
		\item[(P3)] $P$ sends edges in $\mathcal{G}$ to edges or vertices in $\mathcal{F}$. Namely, if $\xymatrix@C=0.5cm{v\ar@{-}[r]&w}$ is an edge, then either $P(v)=P(w)$ or $\xymatrix@C=0.5cm{P(v)\ar@{-}[r]&P(w)}$ is an edge of $\mathcal{F}$;
		\item[(P4)] if $\xymatrix@C=0.5cm{v\ar@{-}[r]&w}$ is an edge in $\mathcal{G}$ such that $v\in \mathcal{F}$ but $w\not\in \mathcal{F}$, then $P(v)=P(w)$.
		
	\end{enumerate}
\end{definition}
The following result is obvious.
\begin{lemma}\label{l:projection-non-leaving}
	Let $\mathcal{G}$ be an exchange graph.
	If there exists a projection for each face $\mathcal{F}$, then $\mathcal{G}$ has non-leaving-face property.
\end{lemma}

 We retain all the notation and conventions of the preceding section. For a cluster $\mathbf{x}_t$, we denote by 
\[[\mathbf{x}_t]=\{x_{1;t},\dots, x_{n,t}\}
\]
the set of cluster variables belonging to $\mathbf{x}_t$, and we refer to $[\mathbf{x}_t]$ the \emph{non-labelled cluster}. We introduce the following relation $R_\mathcal{A}$ on $\mathcal{X}_\mathcal{A}$:
\begin{itemize}
	\item two cluster variables $z_1$ and $z_2$
	have the relation $R_\mathcal{A}$ if they belong to a same non-labelled cluster.
\end{itemize}
According to \cite[Theorem 6.2, Corollary 7.5]{CL2020}, $(\mathcal{X}_\mathcal{A},R_\mathcal{A})$ meets the requirements at the beginning of this section. Hence we have the exchange graph of $(\mathcal{X}_\mathcal{A}, R_\mathcal{A})$.

\begin{definition}
	 The \emph{exchange graph} $\mathcal{G}_\mathcal{A}$ of a cluster algebra $\mathcal{A}$ is defined to be the exchange graph of $(\mathcal{X}_\mathcal{A},R_\mathcal{A})$. Namely, it has non-labelled clusters as vertex set and two non-labelled clusters are joined by an edge if and only if their intersection has cardinality $n-1$.
\end{definition}

\begin{remark}
	\begin{itemize}
		\item[(1)] The above definition is different from \cite[Definition 4.2]{FZ2002} but is equivalent;
		\item[(2)] This note is devoted to studying all cluster algebras and their exchange graphs.  According to \cite[Proposition 6.1]{CL2020}, the exchange graph of a cluster algebra is independent of the choice of its coefficients. Considering the cluster algebra defined in the previous section is sufficient for us.
	\end{itemize}
	
\end{remark}

\subsection{Main result}
In this section we prove that the exchange graph of any cluster algebra has the non-leaving-face property.  Our main idea is to use Bongartz completion to construct the projection needed in Lemma \ref{l:projection-non-leaving}.

\begin{definition}[\cite{CGY2021}]
	Let $U$ be a subset of some non-labelled cluster of $\mathcal{A}$. A non-labelled cluster $B_U(s):=[\mathbf{x}_t]$ is called the \emph{Bongartz completion of $U$ with respect to the vertex $s$} if
	\begin{itemize}
		\item $U$ is a subset of $[\mathbf{x}_t]$;
		\item The $i$-th $c$-vector $\cc_{i;t}^{B_s;s}$ is a non-negative vecotor for any $i$ such that $x_{i;t}\not\in U$.
	\end{itemize}
\end{definition}
The following statement was proved recently in \cite[Theorem 4.14]{CGY2021}.
\begin{theorem}[existence and uniqueness of Bongartz completion]\label{t:bongartz-completion}
Let $U$ be a subset of some non-labelled cluster of $\mathcal{A}$. The Bongartz completion $B_U(s)$ of $U$ with respect to the vertex $s$ exists and is unique.
\end{theorem}

Let's first work with the root vertex $s\in \mathbb{T}_n$.
\begin{lemma}\label{l:case-1}
	If $U\subseteq [\mathbf{x}_s]$, then $B_U(s)=[\mathbf{x}_s]$.
\end{lemma}

\begin{proof}
	Without loss of generality, we may assume that $U=\{x_{1;s}, \dots, x_{l;s}\}$. By definition, it suffices to show that the $i$-th $c$-vector $\cc_{i;s}^{B_s;s}$ is a non-negative vector for any $l+1\leq i\leq n$. But this is obvious, since $C_s^{B_s;s}=I_n$.
\end{proof}

\begin{lemma}\label{l:case-2}
	Let $\xymatrix@C=0.5cm{s\ar@{-}[r]^k&s'}$ be an edge connected to the root vertex $s$.
	 If $U\subseteq [\mathbf{x}_s]$ and $U\not\subseteq [\mathbf{x}_{s'}]$, then $B_U(s')=[\mathbf{x}_s]$.
\end{lemma}
\begin{proof}
	Without loss of generality, we may assume that $U=\{x_{1;s}, \dots, x_{l;s}\}$. By definition, it suffices to show that the $i$-th $c$-vector $\cc_{i;s}^{B_{s'};s'}$ is a non-negative vector for any $l+1\leq i\leq n$. The condition that $U\not\subseteq [\mathbf{x}_{t'}]$ can be stated as
	\[
	1\leq k\leq l<i.
	\]
	Remembering the recursion formula (\ref{f:c-vector}) and $C_{s'}^{B_{s'};s'}=I_n$, we have the following equality 
	\[\cc_{i;s}^{B_{s'};s'}=\cc_{i;s'}^{B_{s'};s'}+[b_{ki;s'}]_+\cc_{k;s'}^{B_{s'};s'}+b_{ki;s'}[-\cc_{k;s'}^{B_{s'};s'}]_+=\ee_i+[b_{ki;s'}]_+\ee_k,
	\]
	which is non-negative. As a consequence, $B_U(s')=[\mathbf{x}_s]$.
\end{proof}

\begin{lemma}\label{l:case-3}
Let \begin{xy}(0,1)*+{t}="A",(10,1)*+{t'}="B",\ar@{-}^k"A";"B" \end{xy} be an edge of $\mathbb{T}_n$. Then either $B_U(t)=B_U(t')$ or $|B_U(t)\cap B_U(t')|=n-1$.
\end{lemma}
\begin{proof}
	Without loss of generality, we may assume that
	$B_U(t)=[\mathbf{x}_v]$ where $v$ is a vertex in $\mathbb{T}_n$ and $U=\{x_{1;v}, \dots, x_{l;v}\}$.
	Let us denote the matrix as $C_v^{B_t;t}=(c_{ij})$. Since $B_U(t)=[\mathbf{x}_v]$, it follows that  $c_{ij}\geq 0$ whenever $l+1\leq j\leq n$. Using Proposition \ref{p:trans-c-vector}, we obtain
	\[C_v^{B_{t'};t'}=(J_k+[-\epsilon_k(C_t^{-B_v;v})B_t]_+^{k\bullet})C_v^{B_t;t}.
	\]More precisely,
	\[
	\begin{aligned}
	C_v^{B_{t'};t'}&=\begin{bmatrix}
	1\\
	& \ddots\\
	[*]_+&\cdots& -1& \cdots & [*]_+\\
	& &  &\ddots\\
	& & &  &1
\end{bmatrix}\begin{bmatrix}
	c_{11}& \cdots& c_{1k}&\cdots& c_{1n}\\
	\vdots& \ddots& \vdots& &\vdots\\
	c_{k1}&\cdots& c_{kk}& \cdots & c_{kn}\\
	\vdots& & \vdots &\ddots&\vdots\\
	c_{n1}&\cdots &c_{nk} &\cdots  &c_{nn}
\end{bmatrix}
\\
&=\begin{bmatrix}
	c_{11}& \cdots& c_{1k}&\cdots& c_{1n}\\
	\vdots& \ddots& \vdots& &\vdots\\
	*&\cdots& *& \cdots & *\\
	\vdots& & \vdots &\ddots&\vdots\\
	c_{n1}&\cdots &c_{nk} &\cdots  &c_{nn}
\end{bmatrix}.
	\end{aligned}
	\]
	
	Let $j\in [l+1,n]$. First we note that if there exists an $i\neq k$ such that $c_{ij}\neq 0$ (hence $c_{ij}>0$), then the sign-coherence of $c$-vectors implies that the $j$-th $c$-vector $\cc_{j;v}^{B_{t'};t'}$ is a non-negative vector. Since $C_v^{B_t;t}$ is invertible over $\mathbb{Z}$, we deduce that there is at most one $j\in [l+1,n]$ such that for any $i\neq k$, $c_{ij}=0$.
	
	If no such $j$ exists, then the last $n-l$ $c$-vectors are all non-negative. Hence \[B_U(t)=[\mathbf{x}_v]=B_U(t').\]
	
	Otherwise, without loss of generality, we may assume that $c_{in}=0$ for all $i\neq k$. In this case, $c_{kn}=1$ for the reason that $C_v^{B_t;t}$ is invertible over $\mathbb{Z}$.
	Now we work with
	 $\xymatrix@C=0.5cm{v\ar@{-}[r]^n&v'}$.
	 Applying (\ref{f:c-vector}), we conclude that the last $n-l$ column vectors of $C_{v'}^{B_{t'};t'}$ are non-negative. That is $[\mathbf{x}_{v'}]=B_U(t')$ and hence $|B_U(t)\cap B_U(t')|=n-1$.
\end{proof}

We are now ready to state and prove the main result of this note.
 \begin{theorem}\label{t:cluster}
	Let $\mathcal{A}$ be a cluster algebra. The exchange graph $\mathcal{G}_\mathcal{A}$ has the non-leaving-face property.
\end{theorem}
\begin{proof}
Let $U$ be a subset of a non-labelled cluster. Let $\mathcal{F}_U$ be the face of $\mathcal{G}_\mathcal{A}$ determined by $U$. Denote by $\text{v}(\mathcal{G}_\mathcal{A})$ (resp. $\text{v}(\mathcal{F}_U)$) the set of vertices of $\mathcal{G}_\mathcal{A}$ (resp. $\mathcal{F}_U$). Theorem \ref{t:bongartz-completion} induces a well-defined map
\begin{eqnarray*}
	P_U:&\text{v}(\mathcal{G}_\mathcal{A})&\to \text{v}(\mathcal{F}_U).\\
	&[\mathbf{x}_t]&\mapsto B_U(t)
\end{eqnarray*}
Lemma \ref{l:case-1}, \ref{l:case-2} and \ref{l:case-3} imply that $P_U$ extends to a projection  $P_U:\mathcal{G}_\mathcal{A}\to \mathcal{F}_U$. Now the result follows from Lemma \ref{l:projection-non-leaving}.
\end{proof}

\begin{remark}
	Instead of using Bongartz completion via $c$-vectors, Theorem \ref{t:cluster} can also be achieved by studying Bongartz co-completion via $g$-vectors \cite{C21}. In particular, the $g$-vector versions of Lemma \ref{l:case-1}, \ref{l:case-2} and \ref{l:case-3} can be obtained from Remark 5.5, Theorem 5.6, Corollary 5.10 and Theorem 5.11 of \cite{C21} along with the $\mathcal{G}$-systems of cluster algebras.
\end{remark}

\noindent{\bf Acknowledgment.}
We thank Dr. Peigen Cao for his valuable comments and discussions.

\bibliographystyle{plain}
\bibliography{NLF-property}

\begin{bibdiv}
\begin{biblist}

\bib{BY2013}{article}{
      author={Br\"ustle, Thomas},
      author={Yang, Dong},
       title={Ordered exchange graphs},
        date={2013},
     journal={Advances in representation theory of algebras},
      volume={EMS Ser. Congr. Rep.,},
       pages={135\ndash 193},
}

\bib{BZ2019}{article}{
      author={Br\"ustle, Thomas},
      author={Zhang, Jie},
       title={Non-leaving-face property for marked surfaces},
        date={2019},
     journal={Front. Math. China},
      volume={14},
      number={3},
       pages={521\ndash 534},
}

\bib{C21}{article}{
      author={Cao, Peigen},
       title={$\mathcal{G}$-systems},
        date={2021},
     journal={Advances in Mathematics},
      volume={377},
      number={107500},
}

\bib{CGY2021}{article}{
      author={Cao, Peigen},
      author={Gyoda, Yasuaki},
      author={Yurikusa, Toshiya},
       title={Bongartz completion via $c$-vectors},
        date={2021},
     journal={preprint},
      volume={arXiv:2106.11668},
}

\bib{CL2020}{article}{
      author={Cao, Peigen},
      author={Li, Fang},
       title={The enough $g$-pairs property and denominator vectors of cluster
  algebras},
        date={2020},
     journal={Mathematische Annalen},
      volume={377},
       pages={1547\ndash 1572},
}

\bib{CP2016}{article}{
      author={Ceballos, Cesar},
      author={Pilaud, Vincent},
       title={The diameter of type d associahedra and the non-leaving-face
  property},
        date={2016},
        ISSN={0195-6698},
     journal={European Journal of Combinatorics},
      volume={51},
       pages={109\ndash 124},
  url={https://www.sciencedirect.com/science/article/pii/S0195669815001055},
}

\bib{FZ2002}{article}{
      author={Fomin, Sergey},
      author={Zelevinsky, Andrei},
       title={Cluster algebras {I: Foundations}},
        date={2002},
     journal={Journal of the American Mathematical Society},
      volume={15},
      number={2},
       pages={497\ndash 529},
}

\bib{FZ2007}{article}{
      author={Fomin, Sergey},
      author={Zelevinsky, Andrei},
       title={Cluster algebras {IV: Coefficients}},
        date={2007},
     journal={Compositio Mathematica},
      volume={143},
      number={1},
       pages={112\ndash 164},
  url={https://www.cambridge.org/core/article/cluster-algebras-iv-coefficients/6BED6C36574A3D8B1297235E607214FE},
}

\bib{GHKK}{article}{
      author={Gross, Mark},
      author={Hacking, Paul},
      author={Keel, Sean},
      author={Kontsevich, Maxim},
       title={Canonical bases for cluster algebras},
        date={2018},
     journal={Journal of the American Mathematical Society},
      volume={31},
       pages={497\ndash 608},
}

\bib{NZ}{article}{
      author={Nakanishi, Tomoki},
      author={Zelevinsky, Andrei},
       title={On tropical dualities in cluster algebras},
        date={2012},
     journal={Algebraic groups and quantum groups},
      volume={Contemp. Math.},
      number={565},
       pages={217\ndash 226},
}

\bib{STT1988}{article}{
      author={Sleator, Daniel~D.},
      author={Tarjan, Robert~E.},
      author={Thurston, William~P.},
       title={Rotation distance, triangulations, and hyperbolic geometry},
        date={1988},
     journal={Journal of the American Mathematical Society},
      volume={1},
      number={3},
       pages={647\ndash 681},
}

\bib{Will2017}{article}{
      author={Williams, Nathan},
       title={W-associahedra have the non-leaving-face property},
        date={2017},
        ISSN={0195-6698},
     journal={European Journal of Combinatorics},
      volume={62},
       pages={272\ndash 285},
  url={https://www.sciencedirect.com/science/article/pii/S0195669817300069},
}

\end{biblist}
\end{bibdiv}

\end{document}